\documentclass[11pt,reqno]{amsbook}
\newtheorem{Thm}{Theorem}[section]

\newtheorem{Lem}[Thm]{Lemma}

\newtheorem{Rmk}[Thm]{Remark}

\usepackage{graphicx}
\usepackage{latexsym}
\begin{document}

\title[On sections of convex bodies] {Modified Busemann-Petty problem on sections  of convex bodies}

\author{A.Koldobsky, V.Yaskin and M.Yaskina}
\address{A.Koldobsky, Department of Mathematics, University of Missouri, Columbia, MO 65211, USA}
\email{koldobsk@math.missouri.edu}

\address{V.Yaskin, Department of Mathematics, University of Missouri, Columbia, MO 65211, USA}
\email{yaskinv@math.missouri.edu}

\address{M.Yaskina, Department of Mathematics, University of Missouri, Columbia, MO 65211, USA}
\email{yaskinam@math.missouri.edu}

\begin{abstract}

The Busemann-Petty problem asks whether origin-symmetric convex bodies in $\mathbb{R}^n$ with smaller central hyperplane sections necessarily
have smaller $n$-dimensional volume. It is known that the answer is affirmative if $n\le 4$ and negative if $n\ge 5$. In this
article we modify the assumptions of the original Busemann-Petty problem to guarantee the affirmative answer in all dimensions.

\end{abstract}

\maketitle

\section{Introduction}

The Busemann-Petty problem asks the following question. Given two  convex origin-symmetric bodies $K$ and $L$ in $\mathbb{R}^n$ such that
$$ \mathrm{vol}_{n-1}(K\cap H)\le \mathrm{vol}_{n-1} (L\cap H)$$
for every central hyperplane $H$ in $\mathbb{R}^n$, does it follow that
$$\mathrm{vol}_n(K)\le \mathrm{vol}_n (L)?$$

The answer to this problem is known to be affirmative if $n\le 4$ and negative if $n\ge 5$.
The solution appeared as the result of work of many mathematicians (see \cite{GKS} or \cite{Z}  for the
solution in all dimensions and historical details).

Since the answer is negative in most dimensions, it is natural to ask what does one need to know about the volumes of central sections  of two
bodies to be able to compare their volumes in all dimensions. Our main result answers this question.

For an origin-symmetric convex body $K$ in $\mathbb{R}^n$, consider the function
$$S_K(\xi)=\mbox{vol}_{n-1}(K\cap\xi^\perp), \qquad \xi\in S^{n-1},$$
where $\xi^\perp $ is the central hyperplane in $\mathbb{R}^n$ orthogonal to $\xi$. We extend $S_K$ from the sphere to the whole $\mathbb{R}^n$
as a homogeneous function of degree $-1$.

Let $\Delta$ be the Laplace operator on $\mathbb{R}^n$. Fractional powers of the Laplacian can be defined as
$$(-\Delta)^{\alpha/2}f =( |x|^\alpha  \hat{f}(x))^\wedge,$$
where  the Fourier transform is considered in the sense of distributions.
If $\alpha$ is an even integer we get the usual Laplacian applied $\alpha/2$ times.

In this article we prove the following
\begin{Thm}\label{Thm:main1}
 Let $\alpha\in (-3,0]$,  $K$ and $L$ be origin-symmetric infinitely smooth convex bodies in $\mathbb{R}^n$,
 $n\ge 4$, so that for every $\xi\in S^{n-1}$
\begin{eqnarray}\label{eqn:condition}
(-\Delta)^{(n-\alpha-4)/2} S_K(\xi)\le (-\Delta)^{(n-\alpha-4)/2} S_L(\xi).
\end{eqnarray}
Then $$\mathrm{vol}_n(K)\le \mathrm{vol}_n(L).$$

On the other hand, for any $\alpha\in(0,1]$ there are convex symmetric bodies $K,L\in
\mathbb{R}^n,$ $n\ge 5$ that satisfy {\rm (\ref{eqn:condition})} for every $\xi\in S^{n-1}$ but
$\mathrm{vol}_n(L)< \mathrm{vol}_n(K).$
\end{Thm}

The result shows that the $\displaystyle\frac{(n-4)}{2}$-th power of the Laplacian is critical
for getting the estimate betweem volumes.

The negative part is formulated only for $\alpha\in (-3,0]$, because we wanted this to work for
$n=5$. In fact, for bigger $n$ one can take bigger values of $\alpha$. Also the condition (1) can
be written in terms of the Fourier transforms so that no smoothness of the bodies is required.

Putting $n=4$ and $\alpha=0$ in the latter theorem one can see that the theorem represents a
generalization of the affirmative part of the solution to the Busemann-Petty problem, and the case
$n=5$ with $\alpha=1$ gives the negative part of the Busemann-Petty problem.

Another generalization of the Busemann-Petty problem was given in \cite{K2}, where the condition
($\ref{eqn:condition}$) was replaced by an inequality for the derivatives of parallel sections
functions at zero. This generalization still involves volumes of non-central sections so it does
not accomplish our goal - to use only central sections to compare volumes. For other
generalizations of the Busemann-Petty problem and related open questions see \cite{BZ}, \cite{K3},
\cite{RZ}, \cite{MP}. In the case where $\alpha=0$ and $n$ is an even integer
the result of Theorem 1.1 was proved in \cite{K4} using an induction argument. The proof from
\cite{K4} can not be extended to other values of $\alpha$ and $n$ and does not produce any results
in the negative direction.

\section{The function $A_{K,\xi,p}$}

Let $K$ be a convex origin-symmetric body in $\mathbb R^n.$  Our definition of a body assumes
that the origin is an interior point of $K$ . The {\it radial function} of $K$ is given by
$$\rho_K(x)=\max \{a>0: ax \in K \}, \ \ \ x\in \mathbb R^n \setminus \{0\} $$

The Minkowski norm of $K$ is defined as $$||x||_K=\min \{a>0: x \in a K \},$$ clearly
$\rho_K(x)=||x||_K^{-1}$.

Writing the volume of $K$ in polar coordinates, one can express the volume in terms of the Minkowski norm:
\begin{eqnarray}\label{polarVolume}
\mathrm{vol}_n(K)=\frac1n \int_{S^{n-1}} ||\theta||^{-n}_K
d\theta.
\end{eqnarray}

We say that a body $K$ is infinitely smooth if its radial function $\rho_K$ belongs to the space $C^\infty(S^{n-1})$ of
infinitely differentiable functions on the unit sphere. Note that  a simple approximation argument reduces the original
Busemann-Petty problem (as well as all generalizations mentioned in the introduction) to the case where the bodies
$K$ and $L$ are infinitely smooth.

Throughout the paper we use the Fourier transform of distributions. The Fourier transform of a distribution $f$ is defined by
$\langle\hat{f}, \hat{\phi}\rangle=(2\pi)^n \langle f, \phi \rangle$ for every test function $\phi$ from the space $ \mathcal{S}$ of rapidly
decreasing infinitely differentiable functions on $\mathbb R^n$.

A distribution is called {\it positive-definite} if for every test function $\phi$
$$\langle f, \phi \ast \overline{\phi(-x)}\rangle \ge 0.$$
By L.Schwartz's generalization of Bochner's theorem, a distribution is positive definite if and only if its Fourier transform is a positive
distribution (in the sense that $\langle \hat{f},\phi \rangle \ge 0$ for every non-negative test function $\phi$; see, for example,
\cite{GV},p.152).

Let $f$ be an integrable continuous function on $\mathbb R$,
$m$-times continuously differentiable in some neighborhood of
zero, $m \in \mathbb{N}.$ For a number  $q\in (m-1,m)$ the {\it
fractional derivative} of the order $q$ of the function $f$ at
zero is defined as follows.

\begin{eqnarray*}
f^{(q)}(0)=\frac{1}{\Gamma(-q)} \int_0^\infty t^{-1-q}\Big( f^{}(t)-f^{}(0)-t f'(0)- \cdots -&& \\
-\frac{t^{m-1}}{(m-1)!}f^{(m-1)}(0)\Big)dt.&&
\end{eqnarray*}

Note that without dividing by $\Gamma(-q)$ the expression for the
fractional derivative represents an analytic function in the
domain $\{q\in\mathbb{C},-1<\mathrm{Re }\ q<m\}$ not including
integers and has simple poles at integers. The function
$\Gamma(-q)$ is analytic in the same domain and also has simple
poles at non-negative integers. Therefore, after division we get
an analytic function on the whole domain
$\{q\in\mathbb{C},-1<\mathrm{Re }\ q<m\}$, which also defines
fractional derivatives of integer orders. Moreover, computing the
limit as $q\to k$, where $k$ is a non-negative integer and $k<m$,
we see that the fractional derivatives of integer orders coincide
with usual derivatives up to a sign:
$$f^{(k)}(0)=(-1)^k\frac{d^k}{dt^k}f(t)|_{t=0}.$$

For $\xi\in S^{n-1},$ consider a function $A_{K,\xi,p}$ on $\mathbb R$ defined by
$$A_{K,\xi,p}(t)=\int_{K\cap\langle x,\xi\rangle=t}|x|^p dx,$$
where $-n+1<p$.

In this section we establish some regularity properties of the
function $A_{K,\xi,p}$ and express its fractional derivatives in
terms of the Fourier transform. We assume $K$ to be an infinitely
differentiable body. In fact this assumption can be weakened if we
require only the existence of finitely many derivatives, as can be
seen from the proof of the following Lemma.

\begin{Lem}\label{Lem:Diff}
Let $\xi \in S^{n-1}$, $k\in \mathbb N$, $-n+k+1<p\le 0.$ Then the function $A_{K,\xi,p}$ is $k$-times continuously differentiable in some
neighborhood of zero. Also, if $q>-1$ is not an integer than the fractional derivative  $A_{K,\xi,p}^{(q)}(0) $  exists if $p>-n+[q]+2$.

If $q\in \mathbb{C}$ then $A_{K,\xi,p}^{(q)}(0) $ is an analytic function of $q$ in the domain $\{q\in \mathbb{C}:\mbox{Re \ }q>-1,
-[\mbox{Re}(-q)]<n+p-1\}$.
\end{Lem}

{\bf Proof.}

First let us prove that $A_{K,\xi,p}(t) $ is continuously differentiable in a neighborhood of zero if $p>-n+2$. Consider a ball of small radius
$s$ centered at zero, that lies entirely in $K$. If $t$ is strictly less than $s$, then the projection of the origin onto the plane $H_t=\{x \in
\mathbb{R}^n : \langle x,\xi\rangle=t\}$ lies inside of the body $K\cap H_t$. Take this point as the origin on the plane $H_t$ and pass to
spherical coordinates. In this coordinate system  we get
$$A_{K,\xi,p}(t)=\int_{S_t^{n-2}}\left( \int_0^{\rho_{K\cap H_t}(\theta)}r^{n-2}(r^2+t^2)^{p/2} dr\right)d\theta $$
where $\rho_{K\cap H_t}(\theta) $ is the radial function of the
body $K\cap H_t$ and $S_t^{n-2}$ is the unit sphere in $H_t$. For
fixed $\xi$ and $\theta,$ we denote by $\rho(t)=\rho_{K\cap
H_t}(\theta) $.

Now let us  show that $\rho(t)$ is differentiable in a
neighborhood of zero and then differentiate with respect to $t$
under the integral. Note that at $t=0$ we can have a problem since
the integral may not converge at 0. Actually, this problem does
not exist if $p=0$, and in this case $A_{K,\xi,p}(t)$ coincides
with the function used in \cite{GKS} which is infinitely
differentiable in a neighborhood of zero.

First let us prove the differentiability of $ \rho(t)$ with
respect to $t$. Consider the two dimensional plane passing through
the origin and spanned by $\theta$ and $\xi$. Let $D$ be the
section of $K$ by this plane, and $\rho_D$ be the radial function
of $D$ defined on $[0,2\pi]$. Since $K$ is infinitely smooth,
$\rho_D$ is $C^\infty$ on the unit circle. Let us get the implicit
formula for $\rho(t)$ in terms of $\rho_D$. The point on the
boundary corresponding to $\rho(t)$ is at the angle $\arctan\left(
\displaystyle\frac{t}{\rho(t)} \right)$ to $\xi$. Therefore, from
the right triangle,
$$\rho(t)=\sqrt{\rho^2_D\left(\arctan\left(\frac{t}{\rho(t)}\right) \right)-t^2 }.$$
By  implicit differentiation $ \rho(t)$ is infinitely differentiable because $K$ is infinitely smooth and contains a neighborhood of the origin
so that the denominator in the formula
$$\rho'(t)=\frac{\rho_D(\arctan(t/\rho))\rho'_D(\arctan(t/\rho))(1/(\rho^2+t^2))\rho-t}{\rho
+\rho_D(\arctan(t/\rho))\rho'_D(\arctan(t/\rho))(1/(\rho^2+t^2))t}$$  is bounded away from 0 uniformly in $\theta$, if $t$ is small enough.

Now that we proved differentiability of $\rho(t)$ let us differentiate under the integral in $\displaystyle\int_0^{\rho(t)}
r^{n-2}(r^2+t^2)^{p/2} dr$. When $t=0$ the integral becomes improper and we need to find the conditions on $p$ to guarantee its convergence.

After differentiation under the integral we get: $$\displaystyle p\int_0^{\rho(t)} r^{n-2}t(r^2+t^2)^{p/2-1} dr.$$

The integrand achieves its maximum in $t$ when  $t^2=\displaystyle\frac{r^2}{1-p}$. Hence the integral above can be estimated as follows
\begin{eqnarray*}
\int_0^{\rho(t)} r^{n-2}t(r^2+t^2)^{p/2-1} dr\le C\int^{\rho(t)}_0r^{n+p-2}dr,
\end{eqnarray*}
the latter being convergent if  $p>-n+2$.

Hence we proved the statement of the Lemma, if $k=1$. Analogously one can show $A_{K,\xi,p}(t)$ is twice continuously differentiable if $p>-n+3$
and, in general, $A_{K,\xi,p}(t)$ is $k$-times continuously differentiable if $p>-n+k+1$.

Now let us consider fractional derivatives of $A_{K,\xi,p}(t)$ at zero. Suppose it is $k$-times continuously differentiable in some small
interval $(-\epsilon,\epsilon)$. If $k-1<q<k$,  its fractional derivative at zero is defined as
\begin{eqnarray*}
A_{K,\xi,p}^{(q)}(0)= \frac{1}{\Gamma(-q)} \int_0^\infty t^{-1-q}\Big(A_{\xi}^{}(t)-A_{\xi}^{}(0)-tA_{\xi}^{'}(0)-\cdots
-&&\\
-\frac{t^{k-1}}{(k-1)!}A_{\xi}^{(k-1)}(0)\Big)dt.&&
\end{eqnarray*}

There may be a problem with convergence of this integral at $t=0$. But in fact if we consider $t\in (0,\epsilon)$, by Taylor's formula,

\begin{eqnarray*}
\int_0^\epsilon t^{-1-q}\left(A_{\xi}^{}(t)-A_{\xi}^{}(0)-tA_{\xi}^{'}(0)-\cdots
-\frac{t^{k-1}}{(k-1)!}A_{\xi}^{(k-1)}(0)\right)dt=&&\\
=\int_0^\epsilon t^{-1-q}\frac{t^{k}}{k!}A_{\xi}^{(k)}(\eta)dt,&&
\end{eqnarray*}
where $\eta\in (0,\epsilon)$ depends on $t$. The latter integral converges since $q<k$.

Now taking $k=[q]+1$ and recalling that $A_{K,\xi,p}(t)\in C^k(-\epsilon,\epsilon)$ if $p>-n+k+1$, one can see that in order to guarantee the
existence of fractional derivatives we need to require that $p>-n+[q]+2$.

Note that both integer and non-integer cases can be written as: $p>-n-[-q]+1$.

If $\{q\in \mathbb{C}:\mbox{Re \ }q>-1, -[\mbox{Re}(-q)]<n+p-1\}$ differentiating the formula for fractional derivatives with respect to $q$,
one can see that $A_{K,\xi,p}^{(q)}(0) $ is an analytic function of $q$.

 \hspace{10.4cm}q.e.d.

The following formula is a generalization of Theorem 2 from
\cite{GKS}.

\begin{Lem}\label{Lem:Main}
 Let $K$ be an infinitely smooth origin-symmetric convex body in $\mathbb{R}^n,$ $q>-1$, $q\ne n+p-1$ and $-n+[q]+2<p \le 0$. Then
 for every $\xi\in S^{n-1},$
$$A_{K,\xi,p}^{(q)}(0)=\frac{\cos\frac{\pi q}{2}}{\pi(n+p-q-1)}(||x||^{-n-p+q+1}\cdot |x|^p_2)^\wedge(\xi).$$

\end{Lem}
{\bf Proof.} Suppose first that $-1<q<0$. The function $\displaystyle A_{K,\xi,p}(t)=\int_{K\cap\langle x,\xi\rangle=t}|x|^p dx=\int_{\langle x,\xi\rangle=t}\chi(||x||)|x|^p dx$ is even. Applying Fubini's theorem and passing to spherical coordinates, we get

\begin{eqnarray*}
A_{K,\xi,p}^{(q)}(0)&=&\frac{1}{\Gamma(-q)}\int_0^\infty t^{-q-1}A_{K,\xi,p}(t)dt\\
                    &=&\frac{1}{2\Gamma(-q)}\int_{-\infty}^\infty |t|^{-q-1}A_{K,\xi,p}(t)dt\\
                    &=&\frac{1}{2\Gamma(-q)}\int_{-\infty}^\infty |t|^{-q-1}\int_{\langle x,\xi\rangle=t}\chi(||x||)|x|^p dx dt\\
                    &=&\frac{1}{2\Gamma(-q)}\int_{\mathbb{R}^n} |\langle x,\xi\rangle|^{-q-1}\chi(||x||)|x|^p dx\\
                    &=&\frac{1}{2\Gamma(-q)}\int_{S^{n-1}} |\langle \theta,\xi\rangle|^{-q-1}\int_0^\infty r^{-q-1}\chi(r||\theta||)r^p r^{n-1} dr
                    d\theta
\end{eqnarray*}
\begin{eqnarray*}
                    &=&\frac{1}{2\Gamma(-q)}\int_{S^{n-1}} |\langle \theta,\xi\rangle|^{-q-1}\int_0^{\frac{1}{||\theta||}}r^{n+p-q-2} dr
                    d\theta\\
                    &=&\frac{1}{2\Gamma(-q)(n+p-q-1)}\int_{S^{n-1}} |\langle \theta,\xi\rangle|^{-q-1}{||\theta||}^{-n-p+q+1} d\theta.
\end{eqnarray*}

Now we extend $A_{K,\xi,p}^{(q)}(0)$ to $\mathbb{R}^n$ as a homogeneous function of $\xi$ of degree $-1-q$. Then for every even test function
$\phi\in \mathcal{S}$,
\begin{eqnarray*}
&&\langle A_{K,\xi,p}^{(q)}(0),\phi(\xi)\rangle  =\frac{1}{2\Gamma(-q)(n+p-q-1)}\times\\
&&\hspace{3cm}\times\int_{S^{n-1}}{||\theta||}^{-n-p+q+1}\int_{\mathbb{R}^n} |\langle \theta,\xi\rangle|^{-q-1} \phi(\xi) d\xi d\theta.
\end{eqnarray*}

Using Lemma 5 from \cite{GKS}
\begin{eqnarray*}
  & = & \frac{-1}{4\Gamma(-q)\Gamma(1+q)(n+p-q-1)\sin\frac{q\pi}{2}}\times\\
  &&\qquad\qquad\qquad\quad \times\int_{S^{n-1}}{||\theta||}^{-n-p+q+1}\int_{-\infty}^\infty |t|^{q}
\hat{\phi}(t\theta) dt d\theta\\
 &=&\frac{-\sin(-\pi q)}{2\pi(n+p-q-1)\sin\frac{q\pi}{2}}\langle(||x||^{-n-p+q+1}\cdot |x|^p_2)^\wedge(\xi),\phi(\xi)\rangle.
\end{eqnarray*}
The latter follows from the fact that $\Gamma(-q)\Gamma(q+1)=-\pi/\sin(q\pi)$ and the calculation
\begin{eqnarray*}
&&\langle(||x||^{-n-p+q+1}\cdot |x|^p_2)^\wedge(\xi),\phi(\xi)\rangle\\
&=&\int_{R^n} ||x||^{-n-p+q+1}\cdot |x|^p_2 \hat{\phi}(x)dx\\
&=&\int_{S^{n-1}}{||\theta||}^{-n-p+q+1}\int_{0}^\infty t^{-n-p+q+1}t^p t^{n-1}\hat{\phi}(t\theta) dt d\theta\\
&=&\int_{S^{n-1}}{||\theta||}^{-n-p+q+1}\int_{0}^\infty t^{q}\hat{\phi}(t\theta) dt d\theta.
\end{eqnarray*}
We have proved that
$$\langle A_{K,\xi,p}^{(q)}(0),\phi(\xi)\rangle=
\frac{\cos\frac{\pi q}{2}}{\pi(n+p-q-1)}\langle(||x||^{-n-p+q+1}\cdot |x|^p_2)^\wedge(\xi),\phi(\xi)\rangle$$ for $-1<q<0$.

To prove the theorem for other values of $q$ we use the fact that  for every even test function $\phi$ the functions
$$ q\mapsto \langle A_{K,\xi,p}^{(q)}(0),\phi(\xi)\rangle $${ and }
$$q\mapsto \frac{\cos\frac{\pi q}{2}}{\pi(n+p-q-1)}\langle(||x||^{-n-p+q+1}\cdot |x|^p_2)^\wedge(\xi),\phi(\xi)\rangle$$ are analytic in the
domain  $\{q\in \mathbb{C}:\mbox{Re\ }q>-1, -[\mbox{Re}(-q)]<n+p-1\}$. The result of the Lemma follows, since these analytic functions coincide
for $q\in (-1,0),$ $\phi$ is arbitrary and, by Lemma \ref{Lem:Diff}, the fractional derivative is a continuous function of $\xi$ outside of the
origin.

\hspace{10.4cm}q.e.d.

\begin{Rmk}  It follows from Lemma \ref{Lem:Main} with $q=\alpha+2$  that $(-\Delta)^{(n-\alpha-4)/2} S_K$ is a real valued function since
up to a coefficient it is equal to $A_{K,\xi,n-\alpha-4}^{(n-\alpha-4)}(0)$. This explains why can
one compare the Laplacians in the statement of Theorem \ref{Thm:main1} .
\end{Rmk}

\begin{Lem}\label{Lem:pos-def}

Let $K$ be an origin-symmetric convex body in $\mathbb{R}^n$. Assume $q\in (-1,2]$ and $-n-[-q]+1<p\le 0$, then $||x||^{-n-p+q+1}\cdot |x|^p_2$
is a positive-definite distribution on $\mathbb{R}^n$.

\end{Lem}

{\bf Proof.} First we prove that
 \begin{eqnarray}\label{eqn:monotA}
 A_{K,\xi,p}^{}(t)\le A_{K,\xi,p}^{}(0), \quad \mbox{ for all } t\ge 0
\end{eqnarray}
If $p=0$ this follows directly from Brunn's theorem (see \cite{S}) stating that the central hyperplane section of a convex body has maximal
volume among all hyperplane sections orthogonal to a given direction. If  $p<0$ one can see that
$$|x|^p=-p\int_0^\infty \chi(z|x|)z^{-p-1}dz,$$ therefore

\begin{eqnarray*}
A_{K,\xi,p}^{}(t)&=&\int_{K\cap\langle x,\xi\rangle=t}|x|^p dx \\
&=&-p\int_{K\cap\langle x,\xi\rangle=t}\int_0^\infty \chi(z|x|)z^{-p-1}dz dx\\
&=&-p\int_0^\infty z^{-p-1}\int_{K\cap\langle x,\xi\rangle=t}\chi(z|x|)dx dz\\
&=& -p\int_0^\infty z^{-p-1}\int_{B_{1/z}\cap K\cap\langle x,\xi\rangle=t} dx dz\\
&\le& -p\int_0^\infty z^{-p-1}\int_{B_{1/z}\cap K\cap\langle x,\xi\rangle=0} dx dz\\
&=&A_{K,\xi,p}^{}(0)
\end{eqnarray*}
by Brunn's theorem applied to the convex origin-symmetric body $B_{1/z}\cap K$, where $B_{1/z}$ is
a ball of radius $\displaystyle\frac{1}{z}$.

Now consider $q\in (1,2)$. Here $\displaystyle \cos\frac{q\pi}{2}$ is negative,  therefore we need
to prove that $A_{K,\xi,p}^{(q)}(0)\le 0$. Using inequality (\ref{eqn:monotA}), the formula for
fractional derivatives for $q\in (1,2)$ and the fact that $A'(0)=0$ we get
\begin{eqnarray*}
A_{K,\xi,p}^{(q)}(0)&=&\frac{1}{\Gamma(-q)}\int_0^\infty t^{-q-1}(A(t)-A(0)-t A'(0))dt \\
&=&\frac{1}{\Gamma(-q)}\int_0^\infty t^{-q-1}(A(t)-A(0))dt\le 0
\end{eqnarray*}
since $\Gamma(-q)$ is positive.

 If $q\in (0,1)$ then $\cos\frac{q\pi}{2}$ is positive and
$$A_{K,\xi,p}^{(q)}(0)=\frac{1}{\Gamma(-q)}\int_0^\infty t^{-q-1}(A(t)-A(0))dt\ge 0$$
since $\Gamma(-q)<0$ for these values of $q$.

Finally if $q\in (-1,0)$ then $\cos\frac{q\pi}{2}$ is positive, $\Gamma(-q)$ is also positive and
$$A_{K,\xi,p}^{(q)}(0)=\frac{1}{\Gamma(-q)}\int_0^\infty t^{-q-1}A(t)dt\ge 0$$

We still have to prove the Lemma for $q= 0, 1, 2.$

When $q=0$, $\cos{\frac{\pi q}{2}} =1$ and
$$A^{(0)}_{K,\xi,p}(0)=(-1)^0 A_{K,\xi,p}(0) \ge 0. $$

When $q=2$, $\cos\frac{\pi q}{2}=-1$ and
$$A^{(2)}_{K,\xi,p}(0)=(-1)^2 A''_{K,\xi,p}(0) \le 0, $$
since $A_{K,\xi,p}(t)$ has maximum at 0.

When $q=1$, take small $\varepsilon >0$. By what we just proved for non-integer $q$, for any non-negative test function $\phi$,
$$\langle(|x|^p_2||x||_K^{-n-p+2+\varepsilon})^\wedge, \phi \rangle\ge 0.$$

Since $||x||_K\le C |x|_2$ for some $C$, it follows that
$$||x||_K^{-n-p+2+\varepsilon}|x|_2^p\le \tilde{C}|x|^{-n+2+\varepsilon} \le \tilde{C}|x|^{-n+1},$$
the latter being a locally-integrable function on $\mathbb{R}^n$.

Set $g(x)=\tilde{C}|x|^{-n+1}|\hat{\phi}(x)|$ for $|x|<1$ and $g(x)=\tilde{C}|\hat{\phi}(x)|$ for $|x|>1$. The function $g(x)$ is integrable on
$\mathbb{R}^n$ and for small $\varepsilon$ we have that $||x||^{-n-p+2+\varepsilon}_K|x|_2^p \hat{\phi}(x)\le g(x)$. Therefore by the Lebesgue
dominated convergence theorem,

$$\langle (||x||^{-n-p+2}_K|x|_2^p )^\wedge,\phi\rangle=\int_{\mathbb{R}^n}||x||^{-n-p+2}_K|x|_2^p \hat{\phi}(x)dx=$$
$$=\lim_{\varepsilon\to 0}\int_{\mathbb{R}^n}||x||^{-n-p+2+\varepsilon}_K|x|_2^p \hat{\phi}(x)dx=
\lim_{\varepsilon\to 0}\langle (||x||^{-n-p+2+\varepsilon}_K|x|_2^p )^\wedge,\phi\rangle\ge 0$$

\hspace{10.4cm}q.e.d.

\section{Proof of Theorem 1.1.}

In this section we prove Theorem \ref{Thm:main1} stated in the introduction.

Let $S_K(\xi)=\mbox{vol}_{n-1}(K\cap\xi^\perp)$, $\xi\in S^{n-1},$
the central section function defined in the Introduction. Then, as
proved in \cite{K1}

\begin{eqnarray}\label{eqn:defS}
S_K(\xi)=\frac{1}{\pi(n-1)}(||x||_K^{-n+1})^\wedge(\xi).
\end{eqnarray}

Extending $S_K(\xi)$ to $\mathbb{R}^n$ as a homogeneous function of degree $-1$ and using the definition of fractional powers of the Laplacian
we get
$$(-\Delta)^{({n-\alpha-4})/{2}}S_L(\theta)= \frac{1}{\pi(n-1)}(|x|^{n-\alpha-4}_2||x||_L^{-n+1})^\wedge(\theta),$$ therefore
\begin{eqnarray*}
&(2\pi)^n &\int_{S^{n-1}}||x||_K^{-1}||x||_L^{-n+1}dx=\\
&=&(2\pi)^n \int_{S^{n-1}}(|x|^{-n+\alpha+4}_2||x||_K^{-1})(|x|^{n-\alpha-4}_2||x||_L^{-n+1})dx\\
&=&\int_{S^{n-1}}(|x|^{-n+\alpha+4}_2||x||_K^{-1})^\wedge(\theta)(|x|^{n-\alpha-4}_2||x||_L^{-n+1})^\wedge(\theta)d\theta\\
&=&\pi(n-1)\int_{S^{n-1}}(|x|^{-n+\alpha+4}_2||x||_K^{-1})^\wedge(\theta)(-\Delta)^{({n-\alpha-4})/{2}}S_L(\theta)d\theta
\end{eqnarray*}
Here we used Parseval's formula on the sphere (see Lemma 3 from
\cite{K2}) and (\ref{eqn:defS}).

By Lemma \ref{Lem:pos-def} with $p=-n+\alpha+4$ and $q=\alpha+2$,
$(|x|^{-n+\alpha+4}_2||x||_K^{-1})^\wedge$ is a non-negative function on $S^{n-1}$, therefore using
the condition of the theorem and repeating the above calculation in the opposite order, we get

$$\int_{S^{n-1}}||x||_K^{-1}||x||_K^{-n+1}dx\le \int_{S^{n-1}}||x||_K^{-1}||x||_L^{-n+1}dx$$

Then by H\"{o}lder's inequality and the polar formula for the
volume (\ref{polarVolume}),

$$n\ \mbox{vol}_n(K)\le \left(\int_{S^{n-1}}||\theta||_K^{-n}d\theta\right)^{1/n} \left(\int_{S^{n-1}}||\theta||_L^{-n}d\theta\right)^{(n-1)/n}=$$
$$n(\mbox{vol}_n(K))^{1/n}(\mbox{vol}_n(L))^{(n-1)/n},$$
which yields the statement of the positive part of the theorem.

Now let us prove the negative part, that is construct two convex
symmetric bodies $K,L\in \mathbb{R}^n,$ $n\ge 5$ such that for
every $\xi$
$$(-\Delta)^{(n-\alpha-4)/2} S_K(\xi)\le (-\Delta)^{(n-\alpha-4)/2} S_L(\xi),$$
but $$\mathrm{vol}_n(L)< \mathrm{vol}_n(K).$$

First assume that $\alpha\in(0,1)$. Again take  $q=\alpha+2$, so $q\in(2,3)$. Let  $p=-n+q+2$. Our
first goal is to construct a body $L$ so that there is a $\xi\in S^{n-1}$ satisfying
\begin{eqnarray}\label{eqn:int}
\int_0^\infty t^{-q-1} \left(A_{L,\xi,p}(t)-A_{L,\xi,p}(0)-A_{L,\xi,p}''(0)\frac{t^2}{2}\right)dt<0.
\end{eqnarray}

Consider the function

$$f(t)=\left(1-t^2-N t^4\right)^\frac{1}{n+p-1}$$

Let $a_N$ be the positive real root of the equation $f(t)=0$. Define the body $L\in \mathbb{R}^n$ as follows.

$$L=\left\{ (x_1,...,x_n)\in \mathbb{R}^n: x_n\in\left[-a_N,a_N\right] \mbox{ and } \left(\sum_{i=1}^{n-1}x_i^2\right)^{1/2}\le f(x_n)\right\},$$
which is a strictly convex infinitely differentiable body.

 Take $\xi$ to be the unit vector in the direction of the $x_n$-axis. Then
 \begin{eqnarray*}
 A_{L,\xi,p}(t)&=&\int_{S^{n-1}}\int_0^{f(t)}(t^2+r^2)^{p/2}r^{n-2}dr\ d\theta\\
 &=&C_n \int_0^{f(t)}(t^2+r^2)^{p/2}r^{n-2}dr
 \end{eqnarray*}
 where $C_n=|S^{n-1}|$.

One can compute:
\begin{eqnarray*}
 A_{L,\xi,p}(0)=\frac{C_n}{n+p-1},
\end{eqnarray*}
and

\begin{eqnarray*}
 A_{L,\xi,p}''(0)=C_n\left[\frac{ p}{n+p-3}-\frac{2}{n+p-1}\right].
\end{eqnarray*}

Now consider those values of $t$ for which $t<f(t)$. Then we can split the integral:
\begin{eqnarray*}
\int_0^{f(t)}(t^2+r^2)^{p/2}r^{n-2}dr=I_1+I_2
\end{eqnarray*}
into two parts, where the first one can be estimated as follows
\begin{eqnarray*}
I_1=\int_0^{t}(t^2+r^2)^{p/2}r^{n-2}dr\le \int_0^{t}(r^2)^{p/2}r^{n-2}dr=\frac{t^{n+p-1}}{n+p-1}
\end{eqnarray*}
and for the second one we will use the inequality:

$$(1+x)^{\gamma}\le 1+\gamma x +\frac{\gamma(\gamma-1)}{2} x^2, \ \ \mbox{ for } \gamma <0 \mbox{ and } 0<x<1.$$

Then

\begin{eqnarray*}
I_2&=&\int_t^{f(t)}(t^2+r^2)^{p/2}r^{n-2}dr\\
&=&\int_t^{f(t)}\left(1+\frac{t^2}{r^2}\right)^{p/2}r^{p+n-2}dr\le\\
&\le& \int_t^{f(t)}\left(1+\frac{p}{2}\frac{t^2}{r^2}+\frac{\frac{p}{2}\left(\frac{p}{2}-1\right)}{2}\frac{t^4}{r^4}\right)r^{p+n-2}dr\\
&=&\left[\frac{r^{n+p-1}}{n+p-1}+\frac{p}{2}\frac{t^2 r^{n+p-3}}{n+p-3}+\frac{\frac{p}{2}\left(\frac{p}{2}-1\right)}{2}\frac{t^4
r^{n+p-5}}{n+p-5}\right]^{f(t)}_t\\
&=&\frac{f^{n+p-1}(t)}{n+p-1}+\frac{p}{2}\frac{t^2 }{n+p-3}f^{n+p-3}(t)+\\
&&+\frac{\frac{p}{2}\left(\frac{p}{2}-1\right)}{2}\frac{t^4 }{n+p-5}f^{n+p-5}(t)
+Ct^{n+p-1}\\
&\le&\frac{f^{n+p-1}(t)}{n+p-1}+\frac{p}{2}\frac{t^2 }{n+p-3}f^{n+p-3}(t)+Ct^{n+p-1}\\
&=&\frac{1-t^2-N t^4}{n+p-1}+\frac{p}{2}\frac{t^2 }{n+p-3}(1-t^2-N t^4)^\frac{n+p-3}{n+p-1}+Ct^{n+p-1}
\end{eqnarray*}
for some constant $C$.

Now we use the inequality:
\begin{eqnarray*}
(1-x)^{\gamma}\ge 1-\gamma x (1-x)^{\gamma-1}, \ \ \mbox{ for } 0<\gamma <1 \mbox{ and } 0<x<1.
\end{eqnarray*}

Therefore,
\begin{eqnarray*}
I_2&\le&\frac{1-t^2-N t^4}{n+p-1}+\frac{p}{2}\frac{t^2 }{n+p-3}\times\\
&\times&\left(1-\frac{n+p-3}{n+p-1}(1-t^2-N t^4)^{\frac{n+p-3}{n+p-1}-1}(t^2+Nt^4)\right)+Ct^{n+p-1}\\
&=&\frac{1-t^2-N t^4}{n+p-1}+\frac{p}{2}\frac{t^2 }{n+p-3}+\\
&&\hspace{2.3cm}+C_1\frac{t^4+Nt^6}{(1-t^2-N t^4)^{\frac{2}{n+p-1}}}+C t^{n+p-1}
\end{eqnarray*}

For the case $t\ge f(t)$ we have the following estimate:

$$\int_0^{f(t)}(t^2+r^2)^{p/2}r^{n-2}dr\le t^p \int_0^{f(t)}r^{n-2}dr=\frac{f^{n-1}(t)}{n-1}t^p\le \frac{t^{n+p-1}}{n-1}$$

In order to estimate the fractional derivative $A^{(q)}_{L,\xi,p}(0)$ we split the integral
(\ref{eqn:int}) into three parts: over $[0,b_N]$, $[b_N, a_N]$ and $[a_N,\infty)$, where $b_N$ is
the positive real root of the equation $1-t^2-N t^4=t^{q+1}$. For the first interval we use the
estimates obtained above for the case $t<f(t)$
\begin{eqnarray*}
&&\int_{0}^{b_N} t^{-q-1} \left(A_{L,\xi,p}(t)-A_{L,\xi,p}(0)-A_{L,\xi,p}''(0)\frac{t^2}{2}\right)dt\le\\
&\le&\int_{0}^{b_N} t^{-q-1} \Big(\frac{1-t^2-N t^4}{n+p-1}+\frac{p}{2}\frac{t^2 }{n+p-3}+C_1\frac{t^4+Nt^6}{(1-t^2-N t^4)^{\frac{2}{n+p-1}}}\\
&&+C
t^{n+p-1}-\frac{1}{n+p-1}-\left[\frac{ p}{n+p-3}-\frac{2}{n+p-1}\right]\frac{t^2}{2}\Big)dt\\
&=&\int_{0}^{b_N} t^{-q-1} \left(\frac{-N t^4}{n+p-1}+C_1\frac{t^4+Nt^6}{(1-t^2-N t^4)^{\frac{2}{n+p-1}}}+C t^{n+p-1}\right)dt
\end{eqnarray*}

Now one can estimate each term of the last integral separately. Since $b_n\approx {N}^{-1/4}$, we get that
$$\int_{0}^{b_N} t^{-q-1} \frac{-N t^4}{n+p-1}dt\approx -C_2 {N}^{q/4}$$
for a positive constant $C_2$.

For the second term, we change the variable of integration: $u=1-t^2-Nt^4$. Then $t\approx\left(\frac{1-u}{N}\right)^{1/4}$,
$dt\approx-N^{-1/4}\frac{du}{4 (1-u)^{3/4}}$. Therefore

\begin{eqnarray*}
&&\int_{0}^{b_N} t^{-q-1}\frac{t^4+Nt^6}{(1-t^2-N t^4)^{\frac{2}{n+p-1}}} dt\\
&&\hspace{4cm}\approx C_3 N^{-1/4}\left(N^{(q-3)/4}+N \cdot N^{(q-5)/4}\right)\\
&&\hspace{4cm}\le C_3 N^{(q-1)/4}
\end{eqnarray*}

 And finally the integral of the last term is small for large values of $N$, since $n+p-1=q+1$. From what we have obtained one can see that the
 integral over $[0,b_N]$ will be negative for large values of $N$ since the leading term is  $-C_2 {N}^{q/4}$.

Now the integral over $[b_N,a_N]$ can be estimated from above by
\begin{eqnarray*}
C\int_{b_N}^{a_N} t^{-q-1}dt\le C\int_{b_N}^{a_N} (b_N)^{-q-1}dt=C \frac{a_N-b_N}{(b_N)^{q+1}}
\end{eqnarray*}
Recalling that $a_N$ and $b_N$ satisfy the equations
$$1-a_N^2-Na_N^4=0 \quad \mbox{ and }\quad 1-b_N^2-N b_N^4=b_N^{q+1}$$ we conclude that
$$b_N^{q+1}=(a_N^2-b_N^2)(1+N(a_N^2+b_N^2)).$$

Therefore
\begin{eqnarray*}
C\int_{b_N}^{a_N} t^{-q-1}dt\le  \frac{C}{(a_N+b_N)(1+N(a_N^2+b_N^2))}\approx C N^{-1/4}.
\end{eqnarray*}

Finally, the integral over $[a_N,\infty)$ can be computed as follows
\begin{eqnarray*}
\int_{a_N}^{\infty} t^{-q-1} \left(-A_{L,\xi,p}(0)-A_{L,\xi,p}''(0)\frac{t^2}{2}\right)dt &\approx& D_1
N^{q/4}+D_2 N^{q/4-1/2}
\end{eqnarray*}
where $D_1<0$. Therefore, this integral is negative for $N$ large enough.

Combining all the  integrals we can see that for $N$ large enough the desired integral (\ref{eqn:int}) is negative. This means that for some
direction $\xi\in S^{n-1}$ the function $(||x||_L^{-1}\cdot |x|^{-n+q+2}_2)^\wedge(\xi)$ is negative, if $q\in (2,3)$.

If $\alpha=1$, so $q=\alpha+2=3$, both sides of the equality in the statement of Lemma
\ref{Lem:Main} vanish, therefore we need to apply the argument from \cite{GKS} (see the proof of
Theorem 1). Then
\begin{eqnarray*}
&&(||x||_L^{-1}\cdot |x|^{-n+5}_2)^\wedge(\xi)=\\
&&\qquad\qquad=C\int_0^\infty t^{-4} \left(A_{L,\xi,p}(t)-A_{L,\xi,p}(0)-A_{L,\xi,p}''(0)\frac{t^2}{2}\right)dt
\end{eqnarray*}
for a positive constant $C$.
Considering the same body as before, we get that $(||x||_L^{-1}\cdot |x|^{-n+5}_2)^\wedge(\xi)$ is also negative at some point $\xi$.

From Lemma \ref{Lem:Main}  the  function $(||x||_L^{-1}\cdot |x|^{-n+q+2}_2)^\wedge(\xi)$ is continuous, hence there is a neighborhood of $\xi$
where it is  negative.

Let $$\Omega=\{\theta \in S^{n-1}: (||x||_L^{-1}\cdot
|x|^{-n+q+2}_2)^\wedge(\theta)<0\}.$$ Choose a non-positive
infinitely-smooth even function $v$ supported on $\Omega$. Extend
$v$ to a homogeneous function $r^{-n+q+1}v(\theta)$ of degree
$-n+q+1$ on $\mathbb{R}^n$. By Lemma 5 from \cite{K2} we know that
the Fourier transform of $ r^{-n+q+1}v(\theta)$ is equal to
$r^{-q-1}g(\theta)$ for some infinitely smooth function $g$ on
$S^{n-1}$.

Define a body $K$ by
$$||x||^{-n+1}_K = ||x||^{-n+1}_L+\varepsilon r^{-n+1} g(\theta)$$
for some small $\varepsilon$ so that the body $K$ is convex. Multiply both sides by
$\displaystyle\frac{|x|^{n-q-2}}{\pi(n-1)}$, apply the Fourier transform and use that $q=\alpha+2$:

\begin{eqnarray*}
(-\Delta)^{(n-\alpha-4)/2}S_K&=&(-\Delta)^{(n-\alpha-4)/2}S_L+\varepsilon\frac{r^{-n+\alpha+3}v(\theta)}{\pi(n-1)} \\
&\le&(-\Delta)^{(n-\alpha-4)/2}S_L,
\end{eqnarray*}
 since $v$ is non-positive.

On the other hand,

\begin{eqnarray*}
\int_{S^{n-1}}(||x||_L^{-1}\cdot |x|^{-n+\alpha+4}_2)^\wedge(\theta) (-\Delta)^{(n-\alpha-4)/2}S_K d\theta=\\
=\int_{S^{n-1}}(||x||_L^{-1}\cdot |x|^{-n+\alpha+4}_2)^\wedge(\theta) (-\Delta)^{(n-\alpha-4)/2}S_L d\theta+\\
+\frac{1}{\pi(n-1)}\int_{S^{n-1}}(||x||_L^{-1}\cdot |x|^{-n+\alpha+4}_2)^\wedge(\theta)\varepsilon r^{-n+\alpha+3}v(\theta) d\theta>\\
>\int_{S^{n-1}}(||x||_L^{-1}\cdot |x|^{-n+\alpha+4}_2)^\wedge(\theta) (-\Delta)^{(n-\alpha-4)/2}S_L
d\theta.
\end{eqnarray*}
Repeating the argument from the proof of the positive part we get:
$$\mathrm{vol}_n(L)<\mathrm{vol}_n(K).$$

 \hspace{10.4cm}q.e.d.

\end{document}